\documentclass[12pt]{amsart}
\usepackage{amssymb}
\usepackage{verbatim}
\usepackage{fullpage}
\theoremstyle{plain}
\newtheorem{theorem}{Theorem}
\newtheorem{corollary}[theorem]{Corollary}
\newtheorem{proposition}[theorem]{Proposition}
\newtheorem{lemma}[theorem]{Lemma}

\theoremstyle{definition}
\newtheorem{definition}[theorem]{Definition}

\theoremstyle{remark}
\newtheorem*{remark}{Remark}

\newcommand{\cstar}{\ensuremath{\text{C}^{*}}-}

\newcommand{\ds}{\displaystyle}
\newcommand{\alg}{\operatorname{Alg}}
\newcommand{\lat}{\operatorname{Lat}}

\newcommand{\man}{\operatorname{Man}}

\newcommand{\supp}{\operatorname{supp}}
\newcommand{\spec}{\operatorname{Spec}}
\newcommand{\spn}{\operatorname{span}}
\renewcommand{\AA}{\mathcal{A}} 
\newcommand{\BB}{\mathcal{B}} 
\newcommand{\sC}{\mathcal{C}} 
\newcommand{\CC}{\mathcal{C}} 
\newcommand{\DD}{\mathcal{D}}
\newcommand{\sH}{\mathcal{H}}
\newcommand{\HH}{\mathcal{H}}
\newcommand{\KK}{\mathcal{K}}
\newcommand{\LL}{\mathcal{L}}
\newcommand{\sM}{\mathcal{M}}
\newcommand{\MM}{\mathcal{M}}

\renewcommand{\SS}{\mathcal{S}}

\newcommand{\bbC}{\mathbb{C}}
\newcommand{\bbN}{\mathbb{N}}

\newcommand{\bbT}{\mathbb{T}}
\newcommand{\bbZ}{\mathbb{Z}}

\newcommand{\fG}{{\mathfrak{G}}}

\newcommand{\fO}{{\mathfrak{O}}}
\newcommand{\fS}{{\mathfrak{S}}}
\newcommand{\fM}{{\mathfrak{M}}}
\newcommand{\bh}{\mbox{${\mathcal B}(\HH)$}}
\newcommand{\cstaralg}{\ensuremath{\text{C}^*}-algebra}
\renewcommand{\H}{{\mathcal{H}}}
\newcommand{\norm}[1]{\left\|{#1}\right\|}
\newcommand{\lip}{\left\langle}
\newcommand{\rip}{\right\rangle}
\newcommand{\eps}{\mbox{$\varepsilon$}}
\renewcommand{\L}{\LL}
\begin{document}

% topmatter
\title[Invariant Manifolds for Operator Algebras]
{Automatic Closure of Invariant Linear Manifolds for Operator Algebras}
 \author{Allan Donsig}
 \address{Department of Mathematics and Statistics\\
        University of Nebraska-Lincoln\\
        Lincoln, NE 68588-0323}
  \email{adonsig@math.unl.edu}
\author{Alan Hopenwasser}
 \address{Department of Mathematics\\
        University of Alabama\\
        Tuscaloosa, AL 35487}
  \email{ahopenwa@euler.math.ua.edu}

\author{David R. Pitts}
   \address{Department of Mathematics and Statistics\\
        University of Nebraska-Lincoln\\
        Lincoln, NE 68588-0323}
  \email{dpitts@math.unl.edu}

 \date{10 May 2000; minor revisions 27 June 2000}
 \keywords{Invariant subspaces, invariant linear manifolds}
 \thanks{2000 {\itshape Mathematics Subject Classification}
	Primary: 47L55; Secondary 47L35, 47L40 }
 \thanks{To appear in {\sl Illinois J. Math.}}

\begin{abstract}
Kadison's transitivity theorem implies that, for irreducible
representations
of \cstar algebras, every invariant linear manifold is closed.
It is known that CSL algebras have this property if, and only if, the 
lattice is hyperatomic (every projection is generated by a finite number 
of atoms).
We show several other conditions are equivalent, including the condition
that every invariant linear manifold is singly generated.

We show that two families of norm closed operator 
algebras have this property.
First, let $\LL$ be a CSL and suppose $\AA$ is a norm closed algebra
which
is weakly dense in $\alg \LL$ and is a bimodule over the 
(not necessarily closed) algebra generated by the atoms of $\LL$.  
If $\LL$ is hyperatomic and the compression of $\AA$ to each atom of
$\LL$
is a \cstar algebra,  then every linear manifold invariant under $\AA$ 
is closed.
Secondly, if $\AA$ is the image of a strongly maximal triangular AF
algebra 
under a multiplicity free nest representation,
where the nest has order type $-\bbN$, then every linear manifold
invariant 
under $\AA$ is closed and is singly generated.
\end{abstract}
\maketitle

The Kadison Transitivity Theorem \cite{rvk57} states, in part, that if
$\pi$ is an irreducible representation of a \cstar algebra $\CC$
acting on a Hilbert space $\HH$, then each linear manifold invariant
under $\pi(\CC)$ is closed.  What other representations also
have the property that every invariant linear manifold is 
closed?  It is not difficult to extend Kadison's result to show that
if $\pi$ is a representation of $\CC$ on $\HH$ then every
invariant linear manifold for $\pi(\CC)$ is closed if, and only if, the
commutant of the image, $\pi(\CC)'$, is a finite dimensional
\cstaralg\ (we outline an argument below).
This condition is in turn equivalent to $\pi(\CC)$ being 
the finite direct sum of irreducible \cstaralg s.  
The summands are unitarily inequivalent if, and
only if, $\pi(\CC)'$ is abelian.  
Thus, if $\pi$ is a multiplicity free representation,
every invariant linear manifold for $\pi(\CC)$ is closed if, and
only if, the lattice of invariant closed subspaces $\lat\pi(\CC)$ is a
finite Boolean algebra.  In the language introduced
below, this says $\lat\pi(\CC)$ is a hyperatomic lattice.

 The main purpose of this note is to give analogous results for
operator algebras which are not \cstar algebras and which have proper
closed invariant subspaces.  Suppose that
$\CC$ is a \cstaralg\ and that $\DD\subseteq\CC$ is a maximal abelian
$*$-subalgebra of $\CC$.  We are interested in representations of
intermediate algebras $\AA$ (so that $\DD\subseteq\AA\subseteq \CC)$.
In this context, a {\it representation of $\AA$} will always be the
restriction of a $*$-representation of $\CC$ to $\AA$.  Such algebras
and their representations have been considered by numerous authors,
including Arveson~\cite{ArvesonSuC*Al,ArvesonSuC*AlII,Arv74}, Muhly,
Qiu and Solel~\cite{MuhQiuSol91}, Peters, Poon and
Wagner~\cite{PetPooWag90}, and Power~\cite{Pow90a}.  In particular,
representations of such algebras have been studied by Orr and
Peters~\cite{OrrPet95} and by Muhly and Solel~\cite{MuhSol96}.
One motivation for considering such representations is that,
under reasonable hypotheses, (general) representations of $\AA$
can often be dilated to $*$-representations of $\CC$; 
see, for example,~\cite{ArvesonSuC*Al,ArvesonSuC*AlII,MS,Don98}.

An important special case is when $\CC=\bh$, $\AA$ is a CSL algebra 
contained in $\CC$, and $\pi$ is the identity representation.
Foia\c{s} \cite{cp71, cp72} determined when every invariant operator 
range for a nest algebra is closed.  
Davidson described invariant operator ranges for reflexive algebras
in~\cite{Dav82}.
These results were extended by the second author to invariant
linear manifolds of CSL algebras in~\cite{ah}, where it is shown that 
every invariant linear manifold for a CSL algebra is closed if, and only if, 
the invariant subspace lattice is hyperatomic.  

The closely related notions of strictly irreducible and topologically
irreducible representations have been studied for Banach algebras; 
see, for example,~\cite{Dix97}.
Also relevant is the transitive algebra problem, which asks if an 
unital operator algebra $\AA$ with $\lat \AA = \{0,I\}$ must
be weakly dense in $\bh$?
An affirmative answer would, of course, also settle the invariant
subspace problem.
It is known that an algebraically transitive subalgebra of $\bh$ must be
weakly dense in $\bh$; see~\cite[Chapter~8]{RadjaviRos73}.
Thus, showing that topological transitivity implies algebraic transitivity
for norm closed operator algebras would also settle the transitive
algebra problem.

Returning to our context, let $\DD\subseteq
\AA\subseteq\CC$ as above and let $\pi$ be a $*$-representation of
$\CC$ such that every invariant linear manifold for $\pi(\AA)$ is
closed.  We wish to observe that  in many situations,
(for example, when $\pi(\DD)''$ is a masa in $\bh$), $\pi(\AA)$ is
$\sigma$-weakly dense in a CSL algebra.  To see this, note that,
\[
\pi(\DD)'' = \alg\lat(\pi (\DD)) \subseteq \alg\lat (\pi(\AA)).
\]
When $\pi(\DD)''$ is a masa in $\bh$ or, more generally,
 $\alg\lat(\pi(\AA))$ contains a masa, then $\alg\lat(\pi(\AA))$ is a
 CSL algebra with invariant subspace lattice $\LL:=\lat(\pi(\AA))$.
 Since every invariant manifold for $\pi(\AA)$ is closed, so is every
 invariant manifold for $\alg\LL$. Therefore, $\LL$ is hyperatomic
 and, in particular, is also atomic.  By~\cite[Theorem~2.2.11]{Arv74},
 $\alg\LL$ is synthetic, and 
then~\cite[Corollary~2 of Theorem~2.1.5]{Arv74}
 shows that $\pi(\AA)$ is $\sigma$-weakly dense in
 $\alg\LL$.

For many of the examples appearing in~\cite{OrrPet95}, $\pi(\DD)''$ is
a masa in $\bh$, so that $\LL:=\lat(\pi(\AA))$ is a commutative
subspace lattice.  Other examples from~\cite{OrrPet95} show that even
when $\pi(\DD)''$ is not a masa, $\alg\lat(\pi(\AA))$ may still
contain a masa, and thus $\alg\lat(\pi(\AA))$ is again a CSL algebra
with lattice $\lat(\pi(\AA))$.

For these reasons, we shall always assume that $\pi(\AA)$ is contained
in the CSL algebra $\alg\LL$, where $\LL=\lat(\pi(\AA))$.  

In section one, we obtain several conditions on a CSL algebra which are
equivalent to the condition that every invariant linear manifold is
closed, and then give an automatic closure result for norm closed
operator algebras which are weakly dense in a CSL-algebra.  In the
second section, we turn to a specific family of norm closed operator
algebras, those arising as representations of triangular AF (TAF)
algebras.

We now turn to a few matters of notation.   All Hilbert spaces in this
paper will be separable.  The symbol $\LL$ always denotes
a CSL, that is, a strongly closed lattice of mutually commuting
projections containing $0$ and $I$.

Given a (not necessarily closed) operator algebra $\AA\subseteq \bh$,
we let $\man \AA$ denote the set of all linear manifolds of $\H$
invariant under $\AA$.  We use $\lat\AA$ for the set of all {\it
closed} subspaces of $\H$ which are invariant under $\AA$.  Clearly
$\lat\AA\subseteq \man\AA$ and the closure of every element of
$\man\AA$ belongs to $\lat\AA$.  Given a vector $x\in\H$, we use
$\{\AA x\}$ for $\{ A x \mid A \in \AA \}$ and $[\AA x]$ for
the closure of $\{\AA x\}$.
For $\AA$ unital, these are the smallest elements of $\man\AA$ 
and $\lat\AA$ respectively containing $x$.

A vector $x\in \H$ is called {\it closed for $\AA$} if $\{\AA x\}=[\AA
x]$. The terminology is justified when one views the vector $x$ as
inducing a map $A\mapsto Ax$ from $\AA$ into $\HH$: a vector is closed
when the associated map has closed range.  (The notion of a closed
vector generalizes the concept of a strictly cyclic vector for an
operator algebra, found in~\cite{FroelichStCyOpAl}: recall that a 
vector $x\in\H$ is strictly cyclic for $\AA$ if $\H=\{\AA x\}$.)  

An element $M\in\lat\AA$ is {\it cyclic} if there exists a vector 
$x\in\H$ such that $M=[\AA x]$. 
To avoid further overloading the word cyclic, we call an invariant
linear manifold $M$ {\it singly generated} if there is a vector 
$x \in\H$ with $M=\{\AA x\}$, whether $M$ is closed or not.

Finally, we outline the argument showing that every invariant
linear manifold for $\pi$, a $*$-representation of  a \cstar algebra
$\sC$,
is closed if, and only if, $\pi(\CC)'$ is finite dimensional.
This is presumably known but we have not found a convenient reference.

One direction is trivial: if $\pi(\CC)'$ is infinite dimensional, it 
contains an infinite chain of projections $P_1 < P_2 < \dots $.  
Then $\sM=\cup_{n=1}^\infty P_n\HH$ is a non-closed invariant manifold 
for $\pi(\CC)$.

For the converse, observe that the finite dimensionality of $\pi(\CC)'$ 
 implies that
that we can decompose $\pi$ as a finite direct sum of irreducible 
representations, say $\oplus \pi_i$, acting on $\oplus \sH_i$.
If $\MM$ is an invariant linear manifold for $\pi(\CC)$, then,
after a possible rearrangement of the order of the summands,  we can
express the elements of $\MM$ as
\[ \bigl(h_1,\ldots,h_k,L_{k+1}(h_1,\ldots,h_k),\ldots,
	L_n(h_1,\ldots,h_k) \bigr) \]
where each $h_i$ is an arbitrary element of $\sH_i$ and
each $L_i$ is a linear transformation from $\oplus_{j=1}^k \sH_j \to
\sH_i$.
If $L$ is the restriction of $L_i$ to some $\sH_j$, $j \le k$, then 
$L$ intertwines the actions of $\pi_j$ and $\pi_i$.

We claim that $L$ is a scalar multiple of a unitary.
Fix a unit vector $x \in \sH_j$.
For each unit vector $v \in \sH_j$, by Kadison's transitivity
theorem, there is a unitary $U \in \sC$ so that $\pi_j(U)x=v$.
As $Lv = L\pi_j(U)x = \pi_i(U) Lx$, we can conclude that
$\|Lv \| = \|Lx\|$ for all unit vectors $v$.
Thus, $\|L\|=\|Lx\|$ for each unit vector $x$ and so
$L$ is a scalar multiple of an isometry.
The transitivity of $\pi_i$ implies that if $L \ne 0$ then $L$ is onto,
so $L$ is a scalar multiple of a unitary, as claimed.
Thus we may write each $L_i$ as a linear combination of unitary
operators.  
It follows that $\sM$ is a closed subspace of $\sH$.

\section{CSL Algebras}

\begin{definition} \label{D:hyper}
A projection $P$ in $\LL$ is \emph{hyperatomic} if there are finitely 
many atoms $A_1, \dots, A_k$ from $\LL$ such that $P = E(A_1, \dots,
A_k)$, 
the smallest projection in $\LL$ containing $A_1, \dots, A_k$.  
We say that $P$ is \emph{generated} by $A_1,\dots, A_k$ if
$P=E(A_1,\dots, A_k)$. 
If each non-zero projection in $\LL$ is hyperatomic, we say that 
the lattice is \emph{hyperatomic}.
\end{definition}

\begin{remark}
A projection $P$ is hyperatomic if, and only if, each ascending
sequence $F_1 \leq F_2 \leq F_3 \leq \dots$ with $P = \vee F_n $, $F_n
\in \LL$, is eventually constant.  Indeed, assume that $P$ is
generated by the atoms $A_1, \dots, A_k $.  Given an increasing
sequence $F_n$ of elements of $\LL$ with $P = \vee F_n $, there is,
for each $j = 1, \dots, n$, a projection $F_{n_j} $ such that $A_j
\leq F_{n_j} $.  If $m \geq
\max\{n_j \mid j = 1, \dots, k\} $ then 
$P =E(A_1, \dots, A_k) \leq F_m \leq P$; hence
$F_m = P$, for all large $m$.

On the other hand, assume the ascending chain condition.  Let $\SS$ be 
the set of atoms contained in $P$.  If $E(\SS) < P $, then
$P - E(\SS) $ contains no atoms, i.e., no non-zero subinterval which is
minimal.  In this case it is easy to construct a strictly increasing
sequence $F_1 < F_2 < \dots$ in $\LL$ such that $P = \vee F_n $,
contradicting the ascending chain condition.  Thus we may assume that
$P$ is generated by the atoms which it contains.  If $P$ is not
generated by finitely many atoms, let $A_1, A_2, \dots$ be an
infinite sequence of atoms which generates $P$.  Then
$E(A_1) \leq E(A_1, A_2) \leq \dots \leq E(A_1, \dots, A_k) < P$, for
all 
$k$, and $P = \vee_k E(A_1, A_2, \dots, A_k)$, again 
contradicting the ascending 
chain condition. Thus we may conclude that $P$ is generated by
finitely many atoms; i.e., $P$ is hyperatomic.  

The version of this remark appropriate to the whole lattice appeared
in \cite{ah}.
\end{remark}

For $x$ and $y$ in $\HH$, $xy^*$ denotes the rank one operator on $\HH$
given by $z\mapsto \lip z,y\rip x$.  
Also, if $P\in\L$, then $P_-$ denotes $\bigvee\{L\in\L: L\not\geq P\}$.

\begin{lemma} \label{L:1atom}
Let $A$ be an atom from $\LL$, let $x$ be a non-zero vector in $A$,
let $y\in E(A)$ and put $T:=\norm{x}^{-2}\, yx^*$.  Then $T\in \alg
\LL$, $Tx = y$ and $TA=T$.  
\end{lemma}

\begin{proof}
 Since $x\in A$, we clearly have $TA=T$.  Notice that if
$L\in\LL$ satisfies $L\not\geq E(A)$, then $AL=0$; for otherwise 
$A\leq L$, whence $E(A)\leq L$.  Therefore
$AE(A)_-=0$, so $x\in (E(A)_-)^\perp$.   Thus, from \cite{wl75}, we see
that $T\in\alg\L$.  
\end{proof}

\begin{proposition} \label{T:xclm}
Let $\LL$ be a commutative subspace lattice  on $\H$ and let
$x\in\H$.    The following are
equivalent:
\begin{enumerate}
\item $x$ is a closed vector for $\alg\LL$, i.e., $\{ \alg\LL x\}$ is
closed.
\item The projection onto $[\alg\LL x]$ is hyperatomic.
\end{enumerate}
\end{proposition}

\begin{proof}
Let $P$ be the orthogonal projection onto $[\alg\LL x]$.  First assume
that $P$ is not hyperatomic.  Let $F_1 < F_2 < \dots$ be a strictly
ascending sequence of projections in $\LL$ such that $P = \vee F_n$.
Let $a_n$ be a sequence of positive real numbers such that
\begin{equation}\label{1.1} \sum_{n=1}^{\infty} n^2 a_n^2 < \infty.
\end{equation}  Let $k_n $ be a
sequence of positive integers such that for all $n\in\bbN$,
\begin{align*}
 &k_{n+1} \geq k_n + 2, \quad \text{and} \\
 &\| ({P - F_{k_n}})x\| = \| F_{k_n}^{\perp}x \| \leq a_n.
\end{align*}
For each $n$, let $y_n \in (F_{k_n + 1} - F_{k_n})\HH $ be a vector with
$\| y_n\| = na_n$.  By (\ref{1.1}), the sum $\sum_{n=1}^{\infty}
y_n$ converges to an element  $y\in P\HH$.

  For every $n\in\bbN$, we have
\[
\| F_{k_n}^{\perp} y\| \geq \|(F_{k_n + 1} - F_{k_n})y \|
   = \|y_n\| = na_n \quad \text{and}\quad 
\|F_{k_n}^{\perp} x \| \leq a_n. 
\]
Hence for all $n$,
\[
\frac {\|F_{k_n}^{\perp}y \|}{\|F_{k_n}^{\perp}x \|} \geq
  \frac {na_n}{a_n} = n.
\]

It follows from \cite{ah80} that $y \neq Tx$ for any $T \in \alg \LL$;
i.e., $y \notin \{\alg\LL x\}$. Thus $\{\alg\LL x\} \neq [\alg\LL x]$,
so $x$ is not a closed vector.

Now suppose that $P$ is hyperatomic.  Let $A_1, \dots, A_n$ be a finite 
set of atoms of $\LL$ such that \[P = E(A_1, \dots, A_n) =
 \bigwedge\{F \in \LL \mid A_j \leq F, j = 1, \dots, n \}.\]  By
deleting
 some atoms, if necessary, we may assume that $A_1, \dots, A_n $ is a
 minimal set which generates $P$.  Thus, if $\SS$ is any proper subset
of 
$\{A_1, \dots, A_n \}$, then $E(\SS) < P $. 

Let $x_k = A_kx$, $k = 1, \dots, n$ .  Note that $x_k \neq 0$, for
each $k$.  
(Otherwise, we have $\overline{\{\alg \LL x\}} \subsetneq P\H=[\alg\LL
x]$
and $x$ is in $E(A_1, \dots, A_{k-1}, A_{k+1}, \dots, A_n)\HH$,
a contradiction.)
Now let $y$ be any vector in $[\alg\LL x]=P\H$.  Then there exist
vectors 
$y_1, \dots, y_n $ (not necessarily unique)
such that $y_k \in E(A_k)\H $, for each $k$, and
$y = y_1 + \dots + y_n$.  
By Lemma \ref{L:1atom}, there is, for each $k$,  an element
 $T_k \in \alg \LL $ such that $T_k x_k = y_k$ and $T_k = T_kA_k$.

Let $T = T_1 + \dots + T_k $.  Clearly, $T_k x_j = 0$ whenever
$k \neq j$.  So $\ds Tx = \sum_{k=1}^{n} T_kx_k = \sum_{k=1}^{n}y_k
 = y$.  This shows that $[\alg\LL x]\subseteq \{\alg\LL x\}$,  and it
follows that $x$ is a closed vector for $\alg\LL$.
\end{proof}

A von Neumann algebra $\fM$ is also a CSL algebra precisely when
$\fM'$ is abelian.  In this sense, one can view CSL
algebras as generalizations of the von Neumann algebras with
abelian commutant.  The discussion in the introduction shows that every
invariant manifold for a von Neumann algebra $\fM$ with abelian
commutant is closed exactly when $\fM$ is finite dimensional, or
equivalently, when $\lat(\fM)$ is hyperatomic.  The next result,
Theorem~\ref{C:singly}, generalizes this characterization to the class
of all CSL algebras.

We also remark that Theorem~\ref{C:singly} extends work of
Froelich in~\cite{FroelichStCyOpAl}.  Motivated by operator theory,
Froelich  introduced the notions of strictly cyclic operator
algebras (those for which there is $x \in \H$ with $\{\AA x\} = \H$)
and of strongly strictly cyclic operator algebras (those for which the
compression of $\AA$ to each invariant projection is strictly cyclic).
He showed that strict cyclicity is equivalent to the ascending chain
condition for the identity and the analogous result for strong strict
cyclicity, essentially $1 \Leftrightarrow 4$ in 
Theorem~\ref{C:singly}.   

\begin{theorem} \label{C:singly}
Let $\LL$ be a commutative subspace lattice.  The following statements
are equivalent.
\begin{enumerate}
\item $\LL$ is a hyperatomic CSL.
\item Every invariant manifold for $\alg\LL$ is closed, i.e., $\man
(\alg\LL)=\LL$.
\item Every singly generated invariant manifold for $\alg\LL$ is closed.
\item Every element of $\LL$ is singly generated; that is, for $P\in\LL$
there exists a vector $x\in P\H$ such that $\{\alg\LL x\}= P\H$.
\item Every invariant manifold for $\alg\LL$ is singly generated.
\end{enumerate}
\end{theorem}

\begin{proof}
($1\Leftrightarrow 2$) This is proved in \cite{ah}.

\begin{comment}
 Let $\MM$ be any invariant linear manifold for
$\alg \LL$ and let $P$ be the projection onto $\overline{\MM}$.  Since
$\LL$ is hyperatomic, we may find a minimal family $A_1,\dots, A_n$ of
atoms of $\LL$ such that such that $P
= E(A_1, \dots, A_n)$.  For each $k$, there is $x_k \in \MM$ such that
$A_k x_k \neq 0$.  Since $A_k x_k \in \MM$, we see that $x :=A_1 x_1 +
\dots + A_nx_n\in\MM$.  Notice that the projection onto $[\alg\LL x]$
is $P$.  Thus by
Proposition~\ref{T:xclm}, $\{\alg\LL x\}=[\alg\LL x] = P\H$.  In
particular, $\MM = [\alg\LL x]$, so $\MM$ is  closed.
  As
$\LL=\lat(\alg\LL)$, we have $\man(\alg\LL)=\LL$.  
\end{comment}

($2\Rightarrow 3$) Obvious.

($3 \Rightarrow 4$) 
 If $P$ is any element of $\LL$ then,
since $\HH$ is separable, there is a vector $x \in \HH$ such that $P$
is the projection onto $[\alg \LL x] $.  But our hypothesis is that
$\{\alg\LL x\}$ is already closed, so $\{\alg\LL x\}=P\H$.

($4\Rightarrow 1$) Given $P\in\LL$ we may find $x\in\HH$ such that
$\{\alg\LL x\}=P\H$, thus $x$ is a closed vector.  By
Proposition~\ref{T:xclm}, $P$ is a hyperatomic projection.  Since $P$
is arbitrary, every projection is hyperatomic and so $\LL$ is
hyperatomic.

($5\Rightarrow 4$) Obvious.

% ($5\Rightarrow 1$) If $\L$ is not hyperatomic, then there exists an
% element $P\in\LL$ which is not a hyperatomic projection.  Let $x\in\H$
% be any vector such that $P\H=[\alg\LL x]$.  By
% Proposition~\ref{T:xclm}, we see that $x$ is not a closed vector, so
% $\{\alg\LL x\}\subsetneq P\H$.  As $x$ was arbitrary, we see that
% $P\H$ is an invariant manifold which is not singly generated.

($2\Rightarrow 5$)  If $\MM$ is any invariant linear manifold, then
$\MM$ is closed by hypothesis, so by the equivalence of (2) and (4),
we see that there exists a vector $x$ such that $\MM=\{\alg\LL x\}$.  
\end{proof}

We next turn our attention to operator algebras $\AA$ which are not
weakly closed, but are subalgebras of CSL algebras. 
Theorem~\ref{C:singly} shows that a necessary condition for $\man(\AA)$
to
coincide with $\lat(\AA)$ is that $\AA$ be a subalgebra contained
inside the algebra of a hyperatomic CSL, so we will restrict our
attention to this setting.

\begin{definition}
Let $\BB\subseteq \bh$ and $\CC\subseteq\bh$ be operator algebras.  
We will say that $\BB$ is {\it $\CC$-transitive} if $\{\BB x\}=\{\CC
x\}$ 
for every $x\in \H$.  
Our primary interest is when $\BB\subseteq \CC$.   
Notice that when $\CC=\bh$, then the statement that $\BB$ is
$\CC$-transitive
is simply the statement that $\BB$ is a transitive operator algebra.  
\end{definition}

The next proposition shows how $\alg\LL$-transitivity, closed vectors,
and automatic closure of invariant manifolds for an algebra
$\AA\subseteq \alg\LL$ are related under the
mild assumption of a ``local approximate unit,'' i.e., when $x\in [\AA
x]$ for every $x\in \HH$. 

\begin{proposition} \label{T:sgen}
Let $\LL$ be a hyperatomic CSL and let $\AA\subseteq \bh$ be an
algebra such that $\lat \AA =\LL$ and such that $x\in [\AA x]$, for
every $x\in\H$.  The following statements are equivalent.
\begin{enumerate}
\item $\AA$ is $\alg\LL$-transitive.
\item Every vector $x\in\H$ is closed for $\AA$.
\item Every invariant manifold for $\AA$ is closed.
\end{enumerate}

Moreover, when any of these conditions hold, every invariant manifold
for $\AA$ is singly generated; i.e., if $\MM\in\man(\AA)$, then there
exists $x\in\HH$ such that $\MM=\{\AA x\}$.
\end{proposition}

Before beginning the proof, we remark that while every element of
$\LL$ is singly generated as an $\alg\LL$ module, it is not {\it a
priori} clear that every element of $\LL$ is singly generated as 
an $\AA$ module.

\begin{proof}
($1\Rightarrow 2$) Let $x\in\H$.  Since $\LL$ is hyperatomic,
Theorem~\ref{C:singly} shows that $[\alg\LL x]=\{\alg\LL x\}$ which,
by assumption, is $\{\AA x\}$.
Thus statement (2) holds.

($2\Rightarrow 1$) By assumption, for all $x$, we have
$\{\AA x\}=[\AA x] \in \LL$. 
Since $\{\alg\LL x\}=[\alg\LL x] $ is the smallest element of
$\LL$ which contains $x$ and $x\in \{\AA x\}$ by hypothesis, we see
$\{\alg\LL x\} \subseteq \{\AA x\}\subseteq \{\alg\LL x\} $.  Thus,
$\{\AA x\} = \{\alg\LL x\} $, for all $x \in \HH$.

($1 \Rightarrow 3$) 
Suppose that $x_1, x_2 \in \HH$.  We claim that there is a vector 
$x$ in $\HH$ such that 
\begin{equation}\label{4_20}
\{\AA x\} = \{\AA x_1\} \vee \{\AA x_2\} 
\end{equation}
Write $P_1$ and $P_2$ for the projections onto 
$\{\alg\LL x_1\} = \{\AA x_1\}$ and $\{\alg\LL x_2\} =  \{\AA x_2\}$.  
These projections are in $\LL$ and hence in $\alg \LL$.  
Let $x = x_1 + P_1^{\perp}x_2$. 
We will show that $\{\alg\LL x\} = \{\alg\LL x_1\} \vee \{\alg\LL
x_2\}$.

Since $x_1 = P_1 x $, it follows immediately that 
$\{\alg\LL x_1\} \subseteq \{\alg\LL x\}$.  
Now let $y \in \{\alg\LL x_2\} \cap \{\alg\LL x_1\}^\perp$.  
Then there is $T \in \alg \LL $ such that $y = Tx_2$.  
Since $x_2 = P_1^{\perp} x_2 + P_1 x_2 $, we have
\[
y = Tx_2 = TP^{\perp}_1 x_2 + TP_1 x_2
 =  TP^{\perp}_1 x_2 + P_1TP_1 x_2 .
\]
Since $P^{\perp}_1 y = y $,
\[
y = P^{\perp}_1 T P^{\perp}_1 x_2 =
  P^{\perp}_1 T P^{\perp}_1 (x_1 + P^{\perp}_1 x_2) =
  P^{\perp}_1 T P^{\perp}_1 x.
\]
This shows that $\{\alg\LL x_2\} \cap \{\alg\LL x_1\}^\perp
 \subseteq \{\alg\LL x\}$.
Combining this with $\{\alg\LL x_1\} \subseteq \{\alg\LL x\}$ gives
$\{\alg\LL x_1\} \vee \{\alg\LL x_2\} \subseteq \{\alg\LL x\}$.  
The reverse inequality follows from the fact that $x = x_1 +
P_1^{\perp}x_2$ 
and the claim is verified.

Now let $\MM$ be an arbitrary invariant linear manifold for $\AA$.  We
need to show that $\MM$ is closed.  Let $Q$ be the projection onto the
closure of $\MM$.  Clearly $Q\in\LL$ and hence is a hyperatomic
projection.  Let $A_1,\dots, A_n$ be a family of atoms of $\LL$ such
that $Q=E(A_1,\dots , A_n)$. Then $Q=E(A_1)\vee E(A_2)\vee \dots\vee
E(A_n)$.  Notice that if $x_j$ is a non-zero vector in $A_j\H$, then
$[\alg\L x_j]=\{\alg\LL x_j\}=\{\AA x_j\}=E(A_j)$.  Inductively applying
(\ref{4_20}), we see that there exists $y\in\H$
such that $\MM=\{\AA y\}$, whence $\MM$ is singly generated.  Since
$\AA$ is $\alg\LL$-transitive, $y$ is a closed vector, whence $\MM$ is
closed.

($3\Rightarrow 2$) Obvious.

It remains to show that when $\man (\AA)=\LL$, then every invariant
manifold for $\AA$ is singly generated.  
Let $\MM$ be an invariant linear manifold for $\AA$.  Then by
hypothesis, the orthogonal projection $Q$ onto $\MM$ belongs to $\LL$,
hence there is a vector $x\in\HH$ such that $\MM=\{\alg\LL x\}$.
Clearly $x\in\MM$. Now $x\in [\AA x]=\{\AA x\}$ and since $\{\alg\LL
x\}$ is the smallest element of $\LL$ containing $x$, we conclude that
$\{\AA x \} \supseteq \{\alg\LL x\}\supseteq \{\AA x\}$.
\end{proof}

The following theorem requires the Kadison transitivity theorem for
its proof and is (partially) an extension of that theorem. 

\begin{theorem} \label{T:dense}
Let $\LL$ be a hyperatomic CSL.  Suppose that $\AA\subseteq \alg\LL$
is a norm closed operator algebra such that
$\overline{\AA}^{\text{wot}}=\alg\LL$.  Assume that $E\AA F \subseteq
\AA$ for all atoms $E$ and $F$ of $\LL$ and that $E\AA E$ is a \cstar
algebra for each atom.  Then $\man(\AA)=\lat(\AA)=\LL$, every
element of $\man(\AA)$ is singly generated, and $\AA$ is
$\alg\LL$-transitive.
\end{theorem}

\begin{proof}
Observe that $E\AA E$ is a \cstaralg\ which is weakly dense in
$E\alg\LL \,E=\BB(E\H)$; thus $E\AA E$ is an irreducible \cstar
-subalgebra of $\BB(E\H)$.

We first assume that the identity operator $I$ is
generated by a single atom $E_0$ of $\LL$.  We shall prove that the
invariant manifold for $\AA$ generated by a unit vector in
$E_0\HH$ is all of $\HH$.  So fix a unit vector $\xi\in E_0\HH$ and
let $x\in\HH$ be any vector.  Let $\{Q_n\}_{n=0}^\infty$ be a sequence
of projections in $\LL''$ such that each $Q_n$ is a finite sum of
atoms of $\LL$, $E_0\leq Q_0$, $\sum_{n=0}^\infty Q_n=I$ and
$\sum_{n=1}^\infty\norm{Q_{n}x} <\infty$.

Fix $n\geq 0$.  Write $Q_n=\sum_{j=1}^{k_n} E_{n,j}$ as a finite sum of  
atoms of $\LL$ and let $x_n=Q_nx$.  Since $E_{n,j}\alg\LL \, 
E_0=\BB(E_0\HH,E_{n,j}\HH)$, and $E_{n,j}\AA E_0$ is weakly dense in  
$E_{n,j}\alg\LL \, E_0$, we may find a norm one 
operator $Y_{n,j}\in\AA$ such  
that $Y_{n,j}=E_{n,j}Y_{n,j}E_0$.  Hence we may find a unit vector  
$u_{n,j}\in E_0\HH$ such that $\norm{Y_{n,j}u_{n,j}}>1/2$.  By  
Kadison's transitivity theorem, there exist unitary operators
$Z_{n,j}\in  
E_{n,j}\AA E_{n,j}$ and $W_{n,j}\in E_0\AA E_0$ such that  
\[\frac{\norm{E_{n,j}x_n}}{\norm{Y_{n,j}u_{n,j}}}Z_{n,j}Y_{n,j}u_{n,j}
=E_{n,j}x_n\quad\hbox{ and }\quad u_{n,j}=W_{n,j}\xi.\] 
 Writing \[A_{n,j}=
\frac{\norm{E_{n,j}x_n}}{\norm{Y_{n,j}u_{n,j}}}Z_{n,j}Y_{n,j}W_{n,j},\]
we see that
$$A_{n,j}\in\AA,\quad \norm{A_{n,j}}<2\norm{E_{n,j}x_n},\quad \hbox{ and
}\quad
  A_{n,j}\xi= E_{n,j}x_n.$$
Therefore, if $B_n=\sum_{j=1}^{k_n} A_{n,j}$, we find $B_n\in \AA$ and  
$B_n\xi =x_n$.  Moreover, since $E_{n,j}A_{n,j}=A_{n,j}$, we find that
for any $\eta\in\HH$,   
$\norm{B_n\eta}^2=\sum_{j=1}^{k_n}\norm{E_{n,j}A_{n,j}\eta}^2$, so
$$\norm{B_n}\leq \left\{ \sum_{j=1}^{k_n} \norm{A_{n,j}}^2  
\right\}^{1/2}<2\norm{x_n}.$$  Notice also that $B_n=B_nE_0$  
by construction.

The fact that $\sum_{n=1}^\infty \norm{x_n} <\infty$ shows that the
series  
$\sum_{n=0}^\infty B_n$ converges uniformly to an element $B\in\AA$.   
Clearly, $B\xi=\sum_{n=0}^\infty B_n\xi=\sum_{n=0}^\infty x_n =x$.
Thus we have shown that the invariant manifold  
generated by $\xi$ is all of $\HH$.   

Furthermore, notice that our construction shows the following:
\begin{enumerate}
\item[a)] $B = BE_0$;
\item[b)] if $E$ is an atom  
of $\LL$ such that $Ex=0$, then $EB=0$; and
\item[c)] $B=\lim_{n\rightarrow \infty} R_n$, where for each $n$,
 $R_n=\sum_{j=1}^{p_n}C_{n,j}$ is a finite sum of elements
$C_{n,j}\in\AA$ which satisfy $C_{n,j} = E_{n,j}C_{n,j}F_{n,j}$ for some
atoms $E_{n,j}$ and
$F_{n,j}$ of $\LL$.
\end{enumerate}

Returning to the general case, if $E$ is any atom from $\LL$, we may
compress to $P(E)$ (i.e., replace $\AA$ by $P(E) \AA P(E)$ acting on
$P(E) \HH$) and apply the argument above to obtain the following:
if $\xi$ is any non-zero vector in $E\HH$ and if $x \in P(E) \HH$,
then there is $B \in \AA$ such that $B\xi = x$, $B = BE$, and $Fx=0$
implies $FB=0$ for all atoms $F$.  (There is one delicate point: our
hypotheses do not guarantee that $P(E) \AA P(E) \subseteq \AA$, but
in the construction of $B$, $B$ is a norm limit of elements which are
finite sums of elements of the form $F_1XF_2$ with $F_2$ and $F_2$
atoms of $\AA$.  Such elements are in $\AA$, by our hypotheses.)

Now let $\MM$ be an invariant linear manifold under $\AA$.
Let $P$ be the  projection onto $\overline{\MM}$.  Then $P$ is invariant
under $\AA$ and, hence, under $\alg \LL$.  So, $P \in \LL$.

Let $E_1, \dots, E_n$ be independent atoms which generate $P$.  
So, $P = P(E_1, \dots E_n)$ and $E_i \alg \LL \,E_j = 0$ whenever
$i \neq j$.  There is a vector $x \in \MM$ such that $E_i x \neq 0$
for all $i$.  (Let $y_i \in E_i \HH$ with $\norm{y_i} = 1$ and
approximate $\sum y_i$ in norm by an element of $\MM$.)

Clearly $\AA x \subseteq \MM$.  We will prove that 
$P\HH \subseteq \AA x$; this implies that $\MM = P\HH$, whence $\MM$
is closed and singly generated.    Let $y \in P\HH$ be arbitrary.  Write
$y = y_1 + \dots + y_n$, where $y_i \in P(E_i)\HH$ for each $i$.
This can be done since $P = \bigvee_i P(E_i)$.

For each $i$, there is an element $B_i \in \AA$ such that
$B_i x_i = y_i$, $B_i = B_i E_i$, and $Fy_i = 0$ implies 
$FB_i =0$, for all atoms $F$.  Let $B = B_1 + \dots + B_n$.
Then
\begin{align*}
Bx &= B_1 x + \dots + B_n x\\
 &= B_1 E_1 x + \dots + B_n E_n x \\
 &= B_1 x_1 + \dots + B_n x_n \\
 &= y_1 + \dots + y_n = y.
\end{align*}
Thus, $y \in \AA \HH$ and $P \HH \subseteq \AA x$.

Finally, since $\AA$ is weakly dense in $\alg\LL$, for every $y\in\HH$
we have $\{\AA y\}=[\AA y]=[\alg \LL y]$.  Since  $\LL$
is hyperatomic, we have $ \{\alg\LL y\}=[\alg\LL y]$, so $\{\AA
y\}=\{\alg\LL y\}$ for every $y\in\H$.  It follows that 
$\AA$ is $\alg\LL$-transitive.
\end{proof}

This theorem  implies immediately a result tacit 
in the proof of the automatic closure theorem in~\cite{ah}:

\begin{corollary}  
Let $\KK$ be the algebra of compact operators and suppose $\LL$ is a 
hyperatomic CSL.  
Then every invariant linear manifold for $\KK \cap \alg \LL$ is closed. 
\end{corollary}

\section{TAF algebras}

We turn now to representations of strongly maximal triangular AF (TAF)
algebras.  These are subalgebras of AF \cstar algebras arising as
limits of triangular digraph algebras and have been extensively
studied; see, for example,~\cite{Pow90a,PetPooWag90,HopPow92,Power92,Don98a}.
If $\AA$ is a closed subalgebra of an AF \cstar algebra $\CC$, then
$\AA$ is \textsl{triangular AF} or \textsl{TAF} if $\AA \cap \AA^*$ 
is a canonical masa in $\CC$.  
A masa $\DD$ is a {\it canonical masa} in $\CC$ if the closed span of 
$N_{\DD}(\CC)$ is $\CC$, where
\[ N_{\DD}(\CC) = \bigl\{ f \in \CC : f \text{ is a partial isometry, }
fdf^*, f^*df \in \DD \text{ for } d \in \DD \bigr\}. \]
A triangular algebra $\AA$ is \textsl{strongly maximal} if
$\overline{\AA+\AA^*}=\CC$.

Let $\AA$ be a strongly maximal triangular AF subalgebra of the AF
\cstaralg\ $\CC$ with $\DD=\AA\cap \AA^*$ a canonical masa in $\CC$.
For reasons we will explain momentarily, we consider representations 
$\pi\colon\AA\rightarrow \bh$ satisfying the following conditions:
\begin{enumerate}
\item  $\pi$ is the restriction to $\AA$ of a $*$-representation $\rho$ 
of $\CC$ on $\HH$; 
\item $\pi(\DD)$ is weakly dense in a masa in $\bh$; and
\item $\lat(\pi(\AA))$ has order type $-\bbN$ and is multiplicity free.
\end{enumerate}
Representations satisfying the first two conditions are called 
\emph{masa preserving}~\cite[p.~130]{OrrPet95}.
Since $\lat(\rho(\CC)) \subseteq \lat(\pi(\AA)) \cap \lat(\pi(\AA^*))
=\{0,I\}$, the $*$-representation $\rho$ is necessarily irreducible.
We will occasionally call a representation which satisfies conditions
1, 2, and 3 an \emph{admissible representation}.

If $\AA$ is non-unital, so is $\CC$.  By $\AA^+$ we mean the obvious
subalgebra of the unitization $\CC^+$ of $\CC$, and it is easy to see
that $\AA^+$ is a strongly maximal TAF subalgebra of $\CC^+$.  Since
$\man\pi(\AA)=\man(\pi(\AA^+)),$ we lose no generality by assuming
that all algebras and representations are unital, and thus we make
this assumption in the sequel. 

The simplest example of a representation satisfying the three conditions
above is the Smith representation of the standard embedding algebra acting on
$\ell^2(-\bbN)$~\cite[Example~I.2]{OrrPet95}.  
In fact, for standard embedding algebras, \cite[Theorem~III.2.1]{OrrPet95} 
shows that $\lat(\pi(\AA))$ is multiplicity free for representations $\pi$
satisfying all the other conditions above.

A more general class of strongly maximal TAF algebras, the
$\bbZ$-analytic algebras considered
in~\cite{OrrPet95,PetPooWag93,PooWag93}, also admit representations of
this form.  However, not all strongly maximal TAF algebras have 
representations satisfying conditions~1, 2, and 3; for example, 
the refinement embedding algebras~(see~\cite{PetPooWag90,OrrPet95}) have 
no such representations.
Further, for a masa preserving representation of a refinement embedding 
algebra, there is a non-closed invariant linear manifold.

We have previously observed that the second condition implies that
$\pi(\AA)$ is $\sigma$-weakly dense in the CSL algebra
$\alg\lat(\pi(\AA))$.  However, since $\AA$ is a strongly maximal TAF
algebra, more is true: \cite[Proposition~0.1]{OrrPet95} shows that the 
second condition implies that $\lat(\pi(\AA))$ is a nest.  
(A nest is a totally ordered CSL.)
Moreover, for many of the examples in \cite{OrrPet95},
$\alg\lat(\pi(\AA))$ is multiplicity free.  If every invariant
manifold for $\pi(\AA)$ is closed, then necessarily the nest
$\lat(\pi(\AA))$ is hyperatomic.  

Furthermore, if $\AA$ is a $\bbZ$-analytic subalgebra of a 
simple AF \cstar-algebra and if $\pi$ is an irreducible 
representation of $C^*(\AA)$ which satisfies condition~2, then 
by~\cite[Proposition~III.3.2]{OrrPet95} $\lat\pi(\AA)$ is a nest 
whose order type is a subset of the integers.  
Since a nest is hyperatomic if, and only if, the complementary nest 
is well-ordered, the hyperatomic nests with order type a subset of the 
integers are just the finite nests and the nests of order type $-\bbN$.  
Automatic closure for invariant manifolds is trivial when $\lat\pi(\AA)$ 
is a multiplicity free finite nest.  
If the nest is finite but not multiplicity free (the nest may be the trivial 
nest $\{0,I\}$, for example), the automatic closure question is open.  
For nests with order type $-\bbN$, Theorem~\ref{TAFrepr} below gives 
an affirmative answer to the question.

Some motivation for, in effect, restricting to irreducible representations
can be found in the following fact, although it does not reduce the study 
of masa preserving representations to the study of irreducible masa preserving 
representations.

\begin{lemma}
Let $\pi$ be a representation of $\CC$ such that $\pi(\DD)$ is weakly dense 
in a masa in $\bh$.  
Every invariant linear manifold for $\pi(\AA)$ is closed if, and only
if, $\pi$ decomposes as a direct sum of finitely many irreducible
representations 
$\pi=\oplus_{i=1}^n\pi_i$ and each invariant linear manifold for
$\pi_i(\AA)$ 
is closed $(i=1,\dots n)$.  
\end{lemma}

\begin{proof}
If $\man(\pi(\AA))=\lat(\pi(\AA))$, then since every invariant manifold 
for $\pi(\CC)$ is also an invariant manifold for $\pi(\AA)$, the
discussion 
in the introduction shows that $\pi$ decomposes as required. 
Then every linear manifold invariant for $\pi_i(\AA)$ is also invariant 
for $\pi(\AA)$.  
Conversely, since $\pi(\DD)''$ is a masa in $\bh$, every invariant
manifold 
$\MM$ for $\pi(\AA)$ decomposes as finite orthogonal sum of invariant
manifolds for $\pi_i(\AA)$, whence $\MM$ is closed.
\end{proof}

We can now state the main theorem of this section.

\begin{theorem}
\label{TAFrepr}
Let $\AA$ be a strongly maximal triangular subalgebra of an
AF \cstar algebra $\CC$.
If $\pi\colon\AA\rightarrow \bh$ is a masa preserving,
order type $-\bbN$, multiplicity free representation, 
then every invariant linear manifold for $\pi(\AA)$ is closed
and singly generated.
\end{theorem}

While the proofs of Theorem~\ref{TAFrepr} and of Theorem~\ref{T:dense} 
employ  similar methods, this theorem is not subsumed by
Theorem~\ref{T:dense} as $\pi(\AA)$ is not a bimodule over the algebra 
generated by the atoms of the nest.  

To prove Theorem~\ref{TAFrepr}, we need to describe admissible 
representations in terms of coordinates.
The full development of such coordinates is technical, and the reader 
is referred to~\cite{MS,Paterson99,Renault80} for more complete 
treatments.
Associated to each AF \cstar algebra $\CC$ there is a unique AF groupoid 
$\fG$ so that the \cstar algebra of $\fG$, $C^*(\fG )$, and
$\CC$ are isomorphic as \cstar algebras.  
The elements of $C^*(\fG)$ can be identified with continuous functions 
on $\fG$.  
With this identification, $C(\fG_0)$ embeds in $C^*(\fG)$ and is a 
canonical masa in $C^*(\fG)$.
In particular, we may identify $\DD$ with $C(\fG_0)$ and $\CC$ with
$C^*(\fG)$.  
Given a unit $e\in \fG_0$, its orbit is the set 
$$[e]:=\{f\in \fG_0: \hbox{ for some } x\in\fG, e=x^{-1}x \hbox{ and }
   xx^{-1}=f\}.$$ 

Given a triangular algebra $\AA$ with $\DD\subseteq \AA\subseteq \CC$,
there is an anti-symmetric subset of $\fG$ containing $\fG_0$, denoted
$\spec(\AA)$, so that $\AA$ is isomorphic to $\{ f \in
C^*(\fG) : \supp f \subset \spec(\AA)\}$. 
If $\AA$ is strongly maximal, then $\spec(\AA)$ totally orders each 
orbit in $\fG_0$.
Similar coordinates can be defined for other groupoid \cstar algebras;
see~\cite{Renault80,Kum86,MuhQiuSol91}.

By Theorem II.1.1 in \cite{OrrPet95}, each representation satisfying
conditions~1 and~2 is unitarily equivalent to a representation of the type
constructed below.
Recall that a {\it 1-cocycle} is a groupoid homomorphism
$\alpha\colon\fG\rightarrow G$, where $G$ is an abelian group; for us
$G$ is $\bbT=\{ x \in \bbC : |x|=1\}$.
As these are the only cocycles we consider, we abbreviate this to
cocycle.

For $v \in N_{\DD}(\CC)$, let $\sigma_v$ denote the partial
homeomorphism
on $\fG_0 = \hat{D}$ induced by the map $d \mapsto vdv^*$.
A measure $\mu$ on $\fG_0$  is \textsl{$\fG$-quasi-invariant}
if, for each $v \in N_{\DD}(\CC)$, the measures $\mu$ and
$\mu\circ\sigma_v$ 
are equivalent, as measures on the domain of $\sigma_v$.
Given such a measure $\mu$, we say a cocycle $\alpha \colon \fG \to \bbT$
is \textsl{$\mu$-measurable} if, for each $v \in N_{\DD}(\CC)$,
the function, denoted $\alpha_v$, from domain of $\sigma_v$ to $\bbC$ 
that sends $x$ to $\alpha(x,\sigma_v(x))$ is measurable.

Since $\CC$ is generated by the diagonal $D \cong C(\fG_0)$ and
$N_D(\CC)$, we can build a representation $\rho$ of $\CC$ on
$L^2(\fG_0,\mu)$ by defining the action of $\rho$ on $\DD$
and on $N_{\DD}(\CC)$ and then extending by linearity to $\CC$.
For $f \in \DD \cong C(\fG_0)$ and $v \in N_{\DD}(\CC)$, define
respectively
\[ \rho(f) \eta = f \eta, \qquad
   \rho(v) \eta = \alpha_v \left[ \frac{ d(\mu \circ \sigma_v)}{d\mu} 
        \right]^{1/2} (\eta \circ \sigma_v). \]

\medskip
\begin{theorem}~\cite[Theorem~II.1.1]{OrrPet95}
Every representation satisfying conditions~1 and~2 is unitarily equivalent
to one arising as above from a $\fG$-quasi-invariant measure $\mu$
and a $\mu$-measurable cocycle $\alpha$, for some choice of $\mu$
and $\alpha$.
\end{theorem}

Suppose now that $\pi$ is an admissible representation.  Since $\pi$
is multiplicity free, the support of $\mu$ is a countable set $S$.
The irreducibility of $\pi$ implies that $S$ is the orbit of a single
point of $\fG_0$ and because $\lat(\pi(\AA))$ has order type $-\bbN$,
$S$ is ordered by $\spec(\AA)$ as $-\bbN$.  Thus $L^2(\fG_0,\mu)$ may
be identified with with $\ell^2(-\bbN)$.  Using $\{e_{j}: j\in
-\bbN\}$ for the basis vectors of $\ell^2(-\bbN)$ and letting $P_n$ be
the projection onto $\overline{\spn\{e_k: k < n\}}$, the lattice of
$\pi(\AA)$ is $\lat\pi(\AA)=\{0,I\}\cup \{P_n:n\in-\bbN\}$.

Given a finite subset $Y \subset \fG_0$, we associate a digraph algebra 
(an algebra isomorphic to a finite-dimensional CSL algebra) to 
$\fS = \spec(\AA) \cap (Y \times Y)$, namely the span of the rank one 
operators $e_x(e_y)^*$ for $(x,y) \in \fS$ acting on the space 
$\ell^2(\{e_y : y \in Y\})$.

\begin{lemma}
\label{covering_alg}
Given a finite subset $Y \subset \fG_0$, let $S$ be the digraph
algebra associated to $\fS = \spec(\AA) \cap (Y \times Y)$.
There is an isometric inclusion $\zeta \colon S \to \AA$ so that
$s$ is in the graph of $\zeta(e_s)$ for each $s \in \fS$.
\end{lemma}

Lemma~\ref{covering_alg} was proved in~\cite[Lemma~4.2]{Don98a};
we need only observe that the inclusion constructed there is 
isometric.

\subsubsection*{Proof of Theorem~\ref{TAFrepr}}
Let $\{P_n: n \in -\bbN \}$ be the projections onto the elements 
of $\lat(\pi(\AA))$, listed in decreasing order; 
thus $0< \ldots < P_{-2}< P_{-1}< P_0 =I$.
For $n\in \{-1,-2,\dots\}$ we let $e_n$ be a unit vector in the range 
of $P_{n+1}-P_n$; since $\lat(\pi(\AA))$ is multiplicity free, 
$\{e_n\}_{n=-\infty}^{-1}$ is an orthonormal basis for $\H$. 

We first show that the singly generated 
invariant manifolds are closed.  Consider the 
manifold $M$ generated by $e_{-1}$.  Clearly, $M$ is dense in $\HH$.
We shall show that if $x\in\HH$ and $\lip x, e_{-1}\rip \neq 0$, then
there exists $T\in\pi(\AA)$ such that $Te_{-1}=x$ and, moreover, that
$T$ can be taken so that $T^{-1}\in \pi(\AA)$.

Fix an element $x\in\HH$  with $\lip x, e_{-1}\rip =1$ and choose a 
decreasing sequence of positive  numbers $\eps_k$ such that 
$\sum_{k=1}^\infty \eps_k =\delta$, where $\delta < (1+\norm{x})^{-1}$.  
Since $P_{-n}^\perp$ are finite rank and converge strongly to $I$, we
may choose an increasing sequence $n_k\in\bbN$  so that 
$\norm{P_{-n_1}x}<1$ and  
$\norm{P_{-n_k}x-P_{-n_{k+1}}x} < \eps_k$, for $k>1$. 

Let $x_1=x-P_{-n_1}x$ and, for $k>1$, let $x_k=
P_{-n_k}x-P_{-n_{k+1}}x$.   
Clearly, $\sum_{k>1}\norm{x_k}<\delta$.
Since there is a natural identification between $(I-P_{-n_k})\HH$
and $\bbC^{n_k}$, we may regard $x_k$ as an element of $\bbC^{n_k}$.

Now let $X_1\in M_{n_1}(\bbC)$ be defined by 
$e_{-1}^{\vphantom{*}}e_{-1}^* + 
 x_1^{\vphantom{*}}e_{-1}^* -I$.  
Here, $e_{-1}$ denotes the ``last standard basis vector'' in
$\bbC^{n_1}$.  
Since $\lip x, e_{-1}\rip =1$, we find that relative to the
decomposition 
$I= (I-e_{-1}^{\vphantom{*}}e_{-1}^*)+e_{-1}^{\vphantom{*}}e_{-1}^*$, 
 $X_1$ has the  upper 
triangular form
$$X_1=\begin{bmatrix} -I_{n_1-1} & v\\ 0& 1\end{bmatrix}.$$ 
Thus $X_1^2=I_{n_1}$.
For $k>1$, let $X_k= x_k^{\vphantom{*}}e_{-1}^* \in T_{n_k}(\bbC)$.

Let $\fO \subset \fG_0$ be the support of the measure $\mu$.
Since $\fO$ has a natural identification with $-\bbN$,
for each $k$, let $Y_k \subset \fO$ be that part identified with
$\{-n_k, \ldots,-2, -1\}$.
Let  $\zeta_n\colon T_{n_k}(\bbC)\rightarrow \AA$ be the 
isometric embedding associated to $Y_k$ given by Lemma~\ref{covering_alg}.
Since  $\zeta_n$ is isometric and 
$\sum_{k>1}\norm{X_k} = \sum_{k>1}\norm{x_k}<\delta $, we see that the
sum
$$\sum_{k=1}^\infty \zeta_{n_k}(X_k)$$ 
converges uniformly to an element $X\in\AA$.  
Notice also that if we let $Z=\sum_{k>1}\zeta_{n_k}(X_k)$, then  
$X=\zeta_{n_1}(X_1) + Z$.  
Since $\zeta_{n_1}(X_1)$ is a square root of $I$, 
$\norm{\zeta_{n_1}(X_1)} < 1+\norm{x}$, and $\norm{Z}<\delta$, we find   
$X=\zeta_{n_1}(X_1)(I + \zeta_{n_1}(X_1) Z)$ is invertible and
$X^{-1}\in\AA$.

Let $T=\pi(X)$.  
Then $T$ is an invertible element of $\pi(\AA)$ and an  
examination of the construction shows that $Te_{-1}= x$.
Note that if $\lip x, e_{-1} \rip = 0$, the same construction still
gives an operator $T$ in $\pi(\AA)$ such that $Te_{-1} = x$;  in this
case $T$ is no longer invertible.

We conclude that if $y_1$ and $y_2$ are vectors in $\HH$ with $\lip
y_1,  
e_{-1}\rip\neq 0$, then there exists $T\in\pi(\AA)$ such that
$Ty_1=y_2$.  
(Indeed,  find   $S_i\in\pi(\AA)$ such that $S_ie_{-1}=y_i$ and $S_1$ is
invertible; then take $T=S_2S_1^{-1}$.)

It follows from our work so far that if $x\in \HH$ has 
$\lip x,e_{-1}\rip\neq 0$, then the invariant manifold 
generated by $x$ is $\HH$, which is obviously closed.

Now let $x\in \HH$ be an arbitrary unit vector and let $M$ be the
invariant manifold generated by $x$.  
The closure of $M$ is an element of the nest, so $\overline{M}=P_{-n}$
for some $n$.  
Clearly, $\lip x, e_{-n}\rip \neq 0$ and by ``compressing'' the argument 
above to $P_{-n}$ we see that  $M=P_{-n}$, so $M$ is closed.  
Hence all singly generated invariant manifolds are closed.

The result now follows from Proposition~\ref{T:sgen}.
\qed

\providecommand{\bysame}{\leavevmode\hbox to3em{\hrulefill}\thinspace}

\end{document}